\documentclass{conm-p-l}

\newtheorem{thm}{Theorem}[section]
\newtheorem{lem}[thm]{Lemma}
\newtheorem{pro}[thm]{Proposition}
\newtheorem{cor}[thm]{Corollary}

\newtheorem{example}[thm]{Example}

\newtheorem{remark}[thm]{Remark}

\newenvironment{enumlist}[1]{%
\begin{list}{}
        {
          
          \settowidth{\labelwidth}{#1}
          \setlength{\leftmargin}{1.1\labelwidth}
        }
      }{%
\end{list}}

\def\R{\mathbb{R}}

\def\T{\mathbb{T}}

\def\2pid{(2 \pi)^{d}}
\def\L2Rd{L^{2}(\R ^{d})}
\def\L2Td{L^{2}(\T ^{d})}

\def\tV{\tilde{V}}
\def\tv{\tilde{v}}
\def\tW{\tilde{W}}

\def\tX{\tilde{X}}
\def\tx{\tilde{x}}
\def\tY{\tilde{Y}}

\def\ty{\tilde{y}}
\def\tz{\tilde{z}}

\def\G{\Gamma}
\def\tG{\widetilde{\Gamma}}
\def\g{\gamma}

\def\Vp{V^{\bot}}
\def\Wp{W^{\bot}}
\def\tVp{\tilde{V}^{\bot}}

\def\sv{\{v_{n}\}}
\def\stv{\{\tv_{n}\}}

\def\la{\langle}
\def\ra{\rangle}
\def\hG{\hat{G}}
\def\hhG{\hat{\hat{G}}}

\def\l2{\ell^{2}}

\begin{document}

\title[Robertson-type theorems]
{Robertson-type Theorems for Countable Groups
of Unitary Operators}

%\makeatletter
\author{David R. Larson}
\address{Department of Mathematics, Texas A\&M University, College Station, TX 77843, USA}
\email{larson@math.tamu.edu}
\thanks{The first and third authors were partially supported by grants from the NSF}

\author{Wai-Shing Tang}
\address{Department of Mathematics, National University of Singapore, 2 Science Drive 2, 117543, Republic of Singapore}
\email{mattws@nus.edu.sg}
\thanks{The second author's research was supported in part by
Academic Research Fund \#R-146-000-073-112, National University of Singapore.}

%    Information for third author
\author{Eric Weber}
\address{Department of Mathematics, Iowa State University, 396 Carver Hall, Ames, IA 50011, USA}
\email{esweber@iastate.edu}
%\thanks{The third author was supported in part by NSF Grant \#000000.}

%    General info
\date{\today }



\copyrightinfo{2006}%            % copyright year
    {American Mathematical Society}% copyright holder

\begin{abstract}
Let $\mathcal{G}$ be a countably infinite group of unitary
operators on a complex separable Hilbert space $H$.
Let $X = \{x_{1},...,x_{r}\}$ and $Y = \{y_{1},...,y_{s}\}$ be
finite subsets of $H$, $r < s$,
$V_{0} = \overline{span} \, \mathcal{G}(X), \, V_1 =
\overline{span} \, \mathcal{G}(Y)$ and
$ V_{0} \subset V_{1} $.
We prove the following result:
Let $W_0$ be a closed linear subspace of $V_1$ such that $V_0 \oplus W_0 = V_1$ (i.e., $V_0 + W_0 = V_1$ and $V_0 \cap W_0 = \{0 \}$).
Suppose that $\mathcal{G}(X)$ and $\mathcal{G}(Y)$
are Riesz bases for $V_{0}$ and $V_{1}$ respectively.
Then there exists a subset $\Gamma =\{z_1,..., z_{s-r}\}$ of $W_0$
such that
$\mathcal{G}(\Gamma)$ is a Riesz basis for $W_0$ if and only if
$ g(W_0) \subseteq W_0 $ for every $g$ in $\mathcal{G}$.
We first handle the case where the group is abelian and then use a cancellation theorem of Dixmier to adapt this 
to the non-abelian case.
Corresponding results for the frame case and the biorthogonal case are also obtained.
\end{abstract}

\maketitle

%%%%%%%%%%%%%%%%%%%%%%%%%%%%%%%%%%%%%%%%%%%%%%%%%%%%%%%%%%%%%%%%%%%%%%%%%%%%%%

\section{Introduction}

\setcounter{equation}{0}

We first consider wavelet-type problems associated with
countably infinite abelian groups of unitary
operators on a complex separable Hilbert space.
Other results in such and similar settings can be found in \cite{DL, HL, HLPS, Weber2}.
We then adapt this to non-abelian groups using a cancellation theorem of Dixmier.

Section~2 is a revisit of Robertson's theorems \cite{Rob} in the
setting of countably infinite abelian groups of unitary
operators.  Recently, Han, Larson,
Papadakis and Stavropoulos \cite{HLPS} extended Robertson's original results to this setting, with proofs involving the spectral
theorem for an abelian group of unitary operators and certain
non-trivial facts on von Neumann algebras.  Using some simple observations
on harmonic analysis,
we show how Robertson's original elementary proofs in \cite{Rob} still
work in this new setting.
In section~3, we
consider the problem of existence of oblique multiwavelets, and show that
oblique multiwavelets exist under a very natural assumption.
The results here extend those of \cite{Tang1, Tang3}.
Finally in section~4, we discuss the non-abelian case using a cancellation theorem of Dixmier.

Let us set up some notations and terminologies.
Throughout this paper, let
$H$ denote a complex separable Hilbert space.  The inner product
of two vectors $x$ and $y$ in $H$ is denoted by $\langle x, y \rangle$.
A countable indexed family ${\sv}_{n \in J}$ of vectors in $H$ is a {\em Riesz basis} for
its closed linear span $V = \overline{span}{\sv}_{n \in J}$
if there exist positive constants $A$ and $B$ such that
\begin{equation}
\label{1.5}
   A \sum |a_{n}|^{2} \leq
   \| \sum a_{n}v_{n} \|^{2} \leq
   B \sum |a_{n}|^{2}, \quad
      \forall \{a_n\} \in \ell^{2}(J).
\end{equation}
A countable indexed family $\sv_{n \in J}$ is a {\em frame} for its closed linear span $V$
if there exist positive constants $A$ and $B$ such that
\begin{equation}
\label{1.4}
  A \|f\|^{2} \leq
  \sum | \langle f, v_{n} \rangle |^{2} \leq B \|f\|^{2}, \quad
  \forall f \in V.
\end{equation}
It is well known that a Riesz basis for a Hilbert space is a frame for the same space.
Two families $\sv$ and
$\stv$ in $H$ are {\em biorthogonal} if
\begin{equation}
\label{1.6}
   \langle v_n, \tv_{m} \rangle = \delta_{n,m} \quad \forall n,\, m\,.
\end{equation}

If $V$ and $W$ are closed linear subspaces of $H$
such that $V \cap W = \{0 \}$  and the vector sum $V_1 = V+W$ is closed,
then we write $V_1 = V \oplus W$ and call this a {\em direct sum}.
In this case, the map $P : V_1 \longrightarrow V_1$ defined by
\begin{equation}
\label{1.7}
  P(v+w) = v, \quad v \in V, \, w \in W,
\end{equation}
is called the (oblique) projection of $V_1$ on $V$ along $W$.
For the special case when $V$ and $W$ are orthogonal, we shall
write $V_1 = V \oplus^{\bot} W$ and call this an {\em orthogonal} direct sum.
We write $V^{\bot}$ for the orthogonal complement of $V$ in $H$.

Let $G$ be a discrete group.
For every $g$ in $G$, let $\chi_g$ denote the
characteristic function of $\{g \}$.  Then $\{\chi_g : g \in G \}$ is an orthonormal basis for $\l2 (G)$.
For each $g$ in $G$, define $l_g : \l2 (G) \longrightarrow \l2 (G)$ by
$(l_{g}a)(h) = a(g^{-1}h), \; h \in G.$
Then $l_g(\chi_h) = \chi_{gh}$ for all $g, \, h$ in $G$.
%For an indexed set $J$,
%the space $\l2 (G \times J)$ can be identified with the spaces
%$\l2 (G) \otimes {\l2}(J)$ and $\l2 (G, \l2 (J))$.
%For each $g$ in $G$, define $L_g : \l2 (G, \l2 (J)) \longrightarrow \l2 (G, \l2 (J))$ by
%$(L_{g}a)(h) = a(g^{-1}h), \; h \in G.$
%Then $L_g = l_g \otimes I$, where $I$ is the identity operator on
%${\l2}(J)$.
The left regular representation $\lambda$ of $G$ is the homomorphism
$\lambda : g \mapsto l_g$.

\section{Robertson's Theorem for Countable Abelian Groups
of Unitary Operators}
\setcounter{equation}{0}

Let $B(H)$ be the space of all bounded linear maps on $H$.
A unitary system $\mathcal{U}$ in $B(H)$ is a set of unitary operators on $H$ which
contains the identity operator on $H$.  A closed linear subspace $M$ of
$H$ is a {\em wandering subspace} for $\mathcal{U}$ if $U(M) \, \bot \, V(M)$
for all $U, V$ in $\mathcal{U}$ with $U \neq V$.
It was first observed in \cite{GLT1} that Robertson's results on
wandering subspaces for the cyclic group generated by a single unitary
operator \cite{Rob} have interesting connection to the existence of orthonormal
wavelets associated with orthonormal multiresolutions.
Subsequently, it was noted in \cite{GLT2} that analogues
of Robertson's results hold for the group generated by a finite set
of commuting unitary operators.  Recently, Han, Larson,
Papadakis and Stavropoulos extended Robertson's results to the
setting of a countable abelian group of unitary operators
\cite[Theorem~4]{HLPS}.  Their proof uses the spectral
theorem for an abelian group of unitary operators and certain
non-trivial facts on von Neumann algebras.  The purpose of this
section is to show, with the help of some simple observations,
how Robertson's original elementary proof in \cite{Rob} can
easily be carried over, almost verbatim, to this new setting.

\begin{lem}
\label{lem:3.1}
Let $M$ and $K$ be wandering subspaces of $H$ for a countable
group $\mathcal{G}$ of unitary operators on $H$.
If $\sum_{g \in \mathcal{G}} \oplus^{\bot} g(M) \subseteq
\sum_{g \in \mathcal{G}} \oplus^{\bot} g(K)$,
then $dim(M) \le dim(K)$.
\end{lem}

\begin{proof}
We modify the arguments due to I. Halperin, given in
\cite[p. 17]{Fill}, for the case of a single unitary operator.
The assertion is trivial if $dim(K) = \infty$.  Therefore suppose that
$dim(K) = k < \infty$.  Let $\{x_1,...,x_k \}$ be an orthonormal basis for
$K$.  Then $\{gx_j : g \in \mathcal{G}, j = 1,...,k \}$ is an orthonormal
basis for $K_1 = \sum_{g \in \mathcal{G}} \oplus^{\bot} g(K)$.
Let $\{y_i \}_{i \in I}$ be an orthonormal basis for $M$, so
$\{gy_i : g \in \mathcal{G}, i \in I \}$ is orthonormal in $K_1$.
By Bessel's inequality,
$$ \sum_{g \in \mathcal{G}} \sum_{i \in I} | \langle x_j, gy_i \rangle |^{2} \;
   \le \| x_j \|^{2} = 1, \quad j = 1,...,k. $$
Hence
\begin{eqnarray*}
 dim(M) & = & \sum_{i \in I} \| y_i \|^{2} \\
        & = & \sum_{i \in I} \sum_{g \in \mathcal{G}}
              \sum_{j=1}^{k} | \langle y_i, gx_j \rangle |^{2} \\
        & = & \sum_{j=1}^{k} \sum_{g \in \mathcal{G}}
              \sum_{i \in I} | \langle x_j, g^{-1}y_i \rangle |^{2} \\
        & \le & k.
\end{eqnarray*}
\end{proof}

Let us recall some results on harmonic analysis (see, for example,
\cite{Rudin}). Let $G$ be a locally compact abelian group, and
 $\hat{G} = \{ \gamma : G \longrightarrow \mathbb{T} \;
  \mbox{continuous characters} \}$
 the dual group of $G$.\\
{\bf Fact 1}: $G$ is isomorphic and homeomorphic to $\hat{\hat{G}}$, the
dual group of $\hat{G}$, via the canonical map $g \longrightarrow e_g$, where
$e_g(\g) = \g (g), \, \g \in \hat{G}$. \\
{\bf Fact 2}: The Fourier transform on $L^1(G)$ is the map
$\wedge : L^1(G) \longrightarrow C_0(\hG)$ defined by
$$\hat{f}(\g) = \int_{G} f(g) \overline{\g (g)} \, dg, \quad \g \in \hG,
\quad f \in L^1(G).$$
Identifying $\hhG$ with $G$ as in Fact~1, the Fourier transform on $L^1(\hG)$ is the map
$\wedge : L^1(\hG) \longrightarrow C_0(\hhG) = C_0(G)$ defined by
$$\hat{f}(g) = \int_{\hG} f(\g) \overline{\g (g)} \, d \g, \quad g \in G,
\quad f \in L^1(\hG).$$\\
{\bf Fact 3}: If $G$ is compact abelian, then $\hat{G}$ forms an orthonormal basis for $L^2(G)$.  Correspondingly, if $G$ is discrete abelian, then
$\hat{\hat{G}} = \{ e_g : g \in G \}$ forms an
orthonormal basis for $L^2(\hat{G})$.\\
{\bf Fact 4}: Suppose that $G$ is a discrete abelian group (so $\hG$ is
compact and $L^2(\hG) \subset L^1(\hG)$).  If $f$ is in $L^2(\hG)$, then
$\hat{f}(g) = \la f, \, e_g \ra , \; g \in G$, and $\hat{f}$ is in
$\ell^{2}(G)$.
Moreover,
\begin{enumlist}{\ (iii) \ }
\item[\ {(i)} \ ]
the Fourier transform $\wedge : L^2(\hG) \longrightarrow \ell^{2}(G)$
is a unitary operator;
\item[\ {(ii)} \ ]
 $\hat{e_g} = \chi_{g}$ for every $g$ in $G$;
\item[\ {(iii)} \ ]
for all $f, h$ in $L^2(\hG)$, we have
$\widehat{fh} = \hat{f} \ast \hat{h}$, i.e.,
$$\widehat{fh}(g) = \sum_{m \in G} \hat{f}(m) \hat{h}(gm^{-1}), \quad g
\in G.$$
\item[\ {(iv)} \ ]
for every $g$ in $G$, we have $\wedge \circ M_g \, = \, l_g \circ \wedge$,
where $M_g : L^2(\hG) \longrightarrow L^2(\hG)$ is the multiplication operator defined by
$$(M_g f)(\g) = \g (g)f(\g), \quad \g \in \hat{G}.$$
\end{enumlist}

Using Fact~$4\,$(i)--(iii), the proofs of Theorem~1 and Theorem~2 in \cite{Rob} can now be
carried over verbatim to our new setting, with $\mathbb{Z}$ and $[0, \, 2 \pi)$
there replaced by $\mathcal{G}$ and $\widehat{\mathcal{G}}$ respectively.
Hence we have

\begin{thm}
\label{thm:3.2}
(\cite[Theorem~4]{HLPS}) Let $\mathcal{G}$ be a countably infinite abelian
 group of unitary operators on $H$,
let $X$ and $Y$ be wandering subspaces of $H$ for $\mathcal{G}$
such that
\begin{enumerate}
\item[(a)]
 $\sum_{g \in \mathcal{G}} \oplus^{\bot} g(X)
\subseteq \sum_{g \in \mathcal{G}} \oplus^{\bot} g(Y)$,
\item[(b)]
 $dim(Y) < \infty$.
\end{enumerate}
Then there exists a wandering subspace $X'$ of $H$ for $\mathcal{G}$ such
that
\begin{enumlist}{\ (iii) \ }
\item[\ {\em (i)} \ ]
$g(X) \perp h(X')$ for all $g, h \in \mathcal{G}$,
\item[\ {\em (ii)} \ ]
$\sum_{g \in \mathcal{G}} \oplus^{\bot} g(X)
\oplus^{\bot}
\sum_{g \in \mathcal{G}} \oplus^{\bot} g(X') =
 \sum_{g \in \mathcal{G}} \oplus^{\bot} g(Y)$,
\item[\ {\em (iii)} \ ]
$dim(X) + dim(X') = dim(Y).$
\end{enumlist}
\end{thm}

Related results for the setting of frames can be found in \cite{HL} and
\cite{Weber2}.

\section{Oblique Multiwavelets and Biorthogonal Multiwavelets}
\setcounter{equation}{0}

Throughout this section,
let $\mathcal{G}$ be a countably infinite abelian group of unitary
operators on a complex separable Hilbert space $H$.
Let $X = \{x_{1},...,x_{r}\}$ and $Y = \{y_{1},...,y_{s}\}$ be
finite subsets of $H$, $r < s$,
$V_{0} = \overline{span} \, \mathcal{G}(X), \, V_1 =
\overline{span} \, \mathcal{G}(Y)$ and
\begin{equation}
\label{4.1}
 V_{0} \subset V_{1}.
\end{equation}

\begin{thm}
\label{thm:4.1}
Let $W_0$ be a closed linear subspace of $V_1$ such that
$V_0 \oplus W_0 = V_1$.
Suppose that $\mathcal{G}(X)$ and $\mathcal{G}(Y)$
are Riesz bases for $V_{0}$ and $V_{1}$ respectively.
Then there exists a subset $\G =\{z_1,..., z_{s-r}\}$ of $W_0$
such that
$\mathcal{G}(\G)$ is a Riesz basis for $W_0$ if and only if
\begin{equation}
  \label{4.20}
g(W_0) \subseteq W_0, \quad g \in \mathcal{G}.
\end{equation}
\end{thm}

Such types of oblique multiwavelets were previously studied in
\cite{Aldrou1, AP, Tang2, Tang3}.
For the proof of Theorem~\ref{thm:4.1}, we need the following
elementary result from \cite[Lemma~3.5]{LTW}.

\begin{lem}
\label{lem:4.3}
Let $M, M'$ and $N$ be linear subspaces of
a vector space $X$ such that
$$X = M \oplus N = M' \oplus N.$$
Let $P$ be the projection of $X$ on $M$ along $N$
and let $Q$ be the projection of $X$ on $M'$ along $N$.
Then $P_1 = P {\vert}_{M'} : M' \longrightarrow M$ and
$Q_1 = Q {\vert}_{M} : M \longrightarrow M'$ are invertible,
and $P_{1}^{-1} = Q_1$.
\end{lem}

\noindent {\em Proof of Theorem~\ref{thm:4.1}}.
The ``only if" part is obvious.  Conversely, suppose that (\ref{4.20}) holds.
Let
$V = V_{1} \cap  V_0^{\bot}$, the orthogonal complement of $V_0$ in $V_1$.
Then we have
$$ V_1 = V_0 \, {\oplus}^{\bot} \, V = V_0 \oplus W_0. $$
By \cite[Proposition~2.2]{LTW} and Theorem~\ref{thm:3.2}, there exists $Z =
\{w_1,..., w_{s-r}\} \subset V$ such that $\mathcal{G}(Z)$ is an
orthonormal basis for $V$.
Let $P : V_1 \longrightarrow V_1$ be the oblique projection of $V_1$
on $W_0$ along $V_0$.
By Lemma~\ref{lem:4.3}, $P {\vert}_{V} : V \longrightarrow W_0$
is invertible.  Hence $\{Pgw_{j} : g \in \mathcal{G} , \, j = 1,...,s-r
\}$ is a Riesz basis for $W_0$.
Since $g(V_0) = V_0$ for every $g \in \mathcal{G}$,
and (\ref{4.20}) holds,
$P$ commutes with $g {\vert}_{V_1}$ for every $g \in \mathcal{G}$.
Therefore,
$\mathcal{G}(\G)$ is a Riesz basis for $W_0$,
where $z_j = Pw_j, \, j = 1,..., s-r$, and
$\G =\{z_1,..., z_{s-r}\}$.
$\: \square$

\vspace{.3cm}

A corresponding result for frames is also valid.
\begin{thm}
\label{thm:4.2}
Let $W_0$ be a closed linear subspace of $V_1$ such that
$V_0 \oplus W_0 = V_1$ (i.e., $V_0 + W_0 = V_1$ and
$V_0 \cap W_0 = \{0 \}$).
Suppose that $\mathcal{G}(X)$ and $\mathcal{G}(Y)$
are frames for $V_{0}$ and $V_{1}$ respectively.
Then there exists a subset $\G =\{z_1,..., z_{p}\}$ of $W_0$
such that
$\mathcal{G}(\G)$ is a frame for $W_0$ if and only if
\begin{equation}
  \label{4.21}
g(W_0) \subseteq W_0, \quad g \in \mathcal{G}.
\end{equation}
\end{thm}

Note that here the number $p$ is not necessarily equal to $s-r$, as in the corresponding case for
Riesz bases.  We omit the proof of Theorem~\ref{thm:4.2}, which is similar to that of Theorem~\ref{thm:4.1}.
Instead of applying Theorem~\ref{thm:3.2}, we use the following simple observation (see \cite{HL}): if 
$\sv_{n \in J}$ is a frame for $V_1$ and $P^{\bot}$ is the orthogonal projection of $V_1$ onto 
$V = V_{1} \cap  V_0^{\bot}$, then $\{ P^{\bot}(v_n) \}_{n \in J}$ is a frame for $V$.

\vspace{.3cm}

As a corollary of the oblique case in Theorem~\ref{thm:4.1}, we obtain
the following result for the biorthogonal case, which is an extension of
\cite[Theorem~3.6]{Tang1}.
\begin{cor}
\label{cor:4.4}
Let $X = \{x_{1},...,x_{r}\},\, \tX = \{\tx_{1},...,\tx_{r}\},\,
Y = \{y_{1},...,y_{s}\}$ and
$\tY = \{\ty_{1},...,\ty_{s}\}$ be
finite subsets of $H$.
Let $\mathcal{G}(X),
\mathcal{G}(\tX),
\mathcal{G}(Y)$ and
$\mathcal{G}(\tY)$ be
Riesz bases for their closed linear spans
$V_{0}, \tV_{0}, V_{1}$ and $ \tV_{1}$ respectively,
$\mathcal{G}(X)$ biorthogonal to $\mathcal{G}(\tX)$,
$\mathcal{G}(Y)$ biorthogonal to $\mathcal{G}(\tY)$, and
$$ V_{0} \subset V_{1}, \quad \tV_{0} \subset \tV_{1}.$$
Let $W_{0} = V_{1} \cap \tVp_0$ and $\tW_0 = \tV_{1} \cap \Vp_{0}.$
If $r < s,$ then
\begin{enumlist}{\ (ii) \ }
\item[\ {\em (i)} \ ]
there exists a subset
$\G :=\{z_1,..., z_{s-r}\}$ of $W_0$
such that $\mathcal{G}(\G)$ is a Riesz basis for $W_0$ and
$\mathcal{G}(X \cup \G)$ is a Riesz basis for $V_1$, and
\item[\ {\em (ii)} \ ]
there exists a subset
$\tG :=\{\tz_1,...,\tz_{s-r}\}$ of $\tW_0$
such that $\mathcal{G}(\tG)$ is a Riesz basis for $\tW_0$ and
$\mathcal{G}(\tX \cup \tG)$ is a Riesz basis for $\tV_1$, and
$\mathcal{G}(\tG)$ is biorthogonal to $\mathcal{G}(\G).$
\end{enumlist}
\end{cor}

\begin{proof}
For every $g \in \mathcal{G}$, we have
$g(V_j) = V_j$ and $g(\tV_{j}) = \tV_{j}$ for $j = 0,1$.
As $g$ is unitary, using the definitions of $W_0$ and $\tW_0$,
$g(W_{0}) = W_{0}$ and $g(\tW_0) = \tW_0$ too.
Since $V_1 = V_0 \oplus W_0$, by Theorem~\ref{thm:4.1},
there exists a subset $\G :=\{z_1,..., z_{s-r}\}$ of $W_0$
such that $\mathcal{G}(\G)$ is a Riesz basis for $W_0$.
By \cite[Proposition~3.1]{Tang1}, $\tW_0 \oplus \Wp_0 = H.$
Hence by \cite[Lemma~3.1]{Tang3},
there exists a subset $\tG :=\{\tz_1,...,\tz_{s-r}\}$ of $\tW_0$
such that $\mathcal{G}(\tG)$ is a Riesz basis for $\tW_0$ and
$\mathcal{G}(\tG)$ is biorthogonal to $\mathcal{G}(\G)$.
The remaining assertions follow easily from \cite[Theoerm~2.1]{Tang1}.
\end{proof}

\section{Non-Abelian Groups}

We consider now the possibility that $\mathcal{G}$ is not abelian.  The proofs of the major statements all follow from the following fundamental result:

\begin{thm}[Cancellation Theorem]
Let $\mathcal{G}$ be a locally compact group, and let $\rho$ be a unitary representation of $\mathcal{G}$ whose commutant is a finite von Neumann algebra.  Suppose $\rho$ is equivalent to $\sigma_{1} \oplus \sigma_{2}$ as well as $\sigma_{1} \oplus \sigma_{3}$.  Then $\sigma_{2}$ is equivalent to $\sigma_{3}$.
\end{thm}

The proof of this theorem is an application of a result due to Dixmier \cite[III.2.3 Proposition 6]{DX}.  We will use the Cancellation Theorem in the following specific case, altered slightly from \cite{BCMO}.

\begin{pro}
Suppose $\mathcal{G}$ is a discrete group, and $\rho$ is a representation of $\mathcal{G}$ which is equivalent to a finite multiple of the left regular representation of $\mathcal{G}$.  Then the commutant of $\rho$ is a finite von Neumann algebra, whence the subrepresentations of $\rho$ satisfy the cancellation property.  That is to say, if $\rho$ is equivalent to $\sigma_{1} \oplus \sigma_{2}$ as well as $\sigma_{1} \oplus \sigma_{3}$, then $\sigma_{2}$ is equivalent to $\sigma_{3}$.
\end{pro}

\begin{thm}
The statements in Theorems 2.2, 3.1, 3.3, and Corollary 3.4 still hold when the assumption that $\mathcal{G}$ is abelian is removed.
\end{thm}

\begin{proof}
We prove Theorem 2.2, without the assumption that $\mathcal{G}$ is abelian.  The other statements follow analogously.

Let $K = \sum_{g \in \mathcal{G}} \oplus^{\perp} g(Y)$ and let $\rho$ denote the action of $\mathcal{G}$ on $K$.  Since $Y$ is a complete wandering subspace of $K$ for $\mathcal{G}$ and is finite dimensional, we have that $\rho$ is equivalent to a finite multiple of the left regular representation $\lambda$ of $\mathcal{G}$, whence we can apply the cancellation theorem.  Let $K_1 = \sum_{g \in \mathcal{G}} \oplus^{\perp} g(X)$; let $\sigma_{1}$ denote the action of $\mathcal{G}$ on $K_1$, and let $\sigma_{2}$ denote the action of $\mathcal{G}$ on $K \cap K_{1}^{\perp}$, the orthogonal complement of $K_1$ in $K$.

For every positive integer $N$, let $\lambda_N$ denote the $N$ multiple of the left regular representation $\lambda$ of $\mathcal{G}$.  We have the following equivalences:
\[ \lambda_{dim(X)} \oplus \lambda_{dim(Y)-dim(X)} \simeq \lambda_{dim(Y)} \simeq \rho \simeq \sigma_{1} \oplus \sigma_{2} \]
as well as
\[ \sigma_{1} \simeq \lambda_{dim(X)} \, , \]
since $X$ is a complete wandering subspace for $\sigma_{1}$ on $K_1$.

Therefore we have that
\[ \sigma_{1} \oplus \lambda_{dim(Y)-dim(X)} \simeq \sigma_{1} \oplus \sigma_{2}. \]
By the cancellation theorem, we have that
\[ \sigma_{2} \simeq \lambda_{dim(Y)-dim(X)} \, , \]
from which it follows that $K \cap K_{1}^{\perp}$ contains a complete wandering subspace of dimension $dim(Y)-dim(X)$; call this subspace $X'$.  The items (i) (ii) and (iii) from Theorem 2.2 now follow.
\end{proof}

%\noindent{\bf Acknowledgement.}
%{\em 
%The second author's research was supported in part by the
%academic Research Fund No. R-146-000-073-112, National University of Singapore.}

\end{document}